\documentclass[a4paper,11pt,fullpage]{article}
\usepackage[english]{babel}

\usepackage[left=1in,right=1in,top=1in,bottom=1in]{geometry}

\usepackage[usenames]{color}
\usepackage{colortbl}

\usepackage{indentfirst}
\usepackage{amsmath}
\usepackage{amsthm}
\usepackage{amsfonts}
\usepackage{amssymb}
\usepackage{amsxtra}
\usepackage{latexsym}
\usepackage{ifthen}
\usepackage{cite}
\usepackage{esint}
\usepackage[cp1250]{inputenc}

\usepackage{epic}
\usepackage{eepic}

\numberwithin{equation}{section}

\def\XXint#1#2#3{{\setbox0=\hbox{$#1{#2#3}{\int}$}
     \vcenter{\hbox{$#2#3$}}\kern-.5\wd0}}

\begin{document}

\title{Qualitative properties of the fourth-order hyperbolic equations}

\author{Kateryna Buryachenko\\
 Humboldt-Universit\"at zu Berlin, Berlin, Germany\\
Vasyl' Stus Donetsk National University, Vinnytsia, Ukraine}



\date{}

\maketitle

\begin{abstract}
We investigate the qualitative properties of the weak solutions to the boundary value problems for the hyperbolic fourth-order linear  equations with constant coefficients in the plane bounded domain convex with respect to characteristics. The main question is to prove the analogue of maximum principle, solvability and uniqueness results for the weak solutions of initial and boundary value problems in the case of weak regularities of initial data from $L^2.$

\end{abstract}


\bigskip
\section{Introduction}
 This paper is devoted to the problem of proving the analog of maximum principle and its further application to the questions of uniqueness and existence for the weak solutions of Goursat, Cauchy and Dirichlet problems for the fourth-order linear hyperbolic equations with the constant coefficients and homogeneous non-degenerate symbol in some plane bounded domain $\Omega\in \mathbb{R}^2$ convex with respect to characteristics:

 \begin{equation}
L(\partial_x)u=a_0\frac{\partial^4u}{\partial x_1^4}+a_1\frac{\partial^4u}{\partial x_1^3\partial x_2}+a_2\frac{\partial^4u}{\partial x_1^2\partial x_2^2}+a_3\frac{\partial^4u}{\partial x_1\partial x_2^3}+a_4\frac{\partial^4u}{\partial x_2^4}=f(x).\label{eq:01}
\end{equation}
 Here coefficients $a_j,\,j=0,\,1,...,\,4$ are constant, $f(x)\in L^2(\Omega),$ $\partial_x=\left(\frac{\partial}{\partial x_1},\,\frac{\partial}{\partial x_2}\right).$ We assume, that Eq. (\ref{eq:01}) is hyperbolic, that means that all roots of characteristics equation
 \begin{equation*}
 L(1,\,\lambda)=a_0\lambda^4+a_1\lambda^3+a_2\lambda^2+a_3\lambda +a_4=0
 \end{equation*}
 are  prime, real  and are not equal to  $\pm i$, that means that the symbol of Eq. (\ref{eq:01}) is non-degenerate or that the Eq. (\ref{eq:01}) is a principal-type equation. The equations for which the roots of the corresponding characteristic equation are multiple and can take the values $\pm i$ are called the equation with
degenerate symbol (see \cite{Bur1}).

 The main novelty of the paper is to prove the analog of maximum principle for the fourth-order hyperbolic equations. This question is very important due to usually a natural physical interpretation, and that it helps to establish the qualitative properties of the solutions (in our case the questions of uniqueness and existence of weak solution). But as it is well known, the maximum principle even for the simple case of hyperbolic equation (one dimensional wave equation \cite{Protter}) are quite different from those for elliptic and parabolic cases,  for which it is a natural fact, such a way a role of characteristics curves and surfaces becomes evident in  the situation of hyperbolic type PDEs.

 We call the angle of characteristics slop the solution to the equation $-\tan\varphi_j=\lambda_j$, and angle between $j-$ and $k-$ characteristics: $\varphi_k-\varphi_j\neq \pi  l,\, l\in \mathbb{Z},$     where $\lambda_j\neq\pm i$ are real and prime roots of the characteristics equation, $j,\,k=1,\,2,\,3,\,4$.

 Most of these equations serve as  mathematical models of many physical processes and attract the interest of researchers. The most famous of them are elasticity beam equations (Timoshenko beam equations with and without internal damping) \cite{Capsoni}, short laser pulse equation \cite{Fichera}, equations which describe the structures are subjected to moving loads, and equation of Euler-Bernoulli beam resting on two-parameter Pasternak foundation and subjected to a moving load or mass \cite{Uzz} and others.

Due to evident practice application, these models need more exact tools for studying, and as consequence, to attract fundamental knowledge. As usual, most of these models are studied by the analytical-numerical methods (Galerkin's methods).

The range of problems studied in this work belongs to a class of quite actual problems of well-posedness of so-called general boundary-value problems for higher-order differential equations originating from the works by L. Hormander and M.Vishik who used the theory of extensions to prove the existence of well-posed boundary-value problems for linear differential equations of arbitrary order with constant complex coefficients in a bounded domain with smooth boundary. This theory got its present-day development in the works by G. Grubb \cite{Grubb}, L.Hormander \cite{H1}, and A. Posilicano \cite{Pos}. Later, the problem of well-posedness of boundary-value problems for various types of second order differential equations was studied by V. Burskii and A. Zhedanov \cite{Bursk1}, \cite{Bursk2} which developed a method of traces associated with a differential operator and applied this method to establish the Poncelet, Abel and Goursat problems.
 In the previous works of author (see \cite{Bur3}) there have been developed qualitative methods of studying  Cauchy problems and nonstandard in the case of hyperbolic equations Dirichlet and Neumann problems for the linear fourth-order equations (moreover, for an equation  of any even order $2m,\, m\geq 2,$ ) with the help of operator methods (L-traces, theory of extension, moment problem, method of duality equation-domain and others), \cite{Bur2}.
 There were proved the existence and uniqueness results,  obtained the criteria of nontrivial solvability of the Dirichlet and Neumann problems in a disk for the principal type equations and equations with degenerate symbol, in particular, there were established the interrelation
between the multiplicity of roots of characteristics equation and the existence of a nontrivial solution of the corresponding Dirichlet and Neumann problems (as a consequences, there was established the Fredholm property of the operator
of such problems).

As concern maximum principle, at the present time there are not any results for the fourth order equations even in linear case. As it was mentioned above, the maximum principle even for the simple case of one dimensional wave equation \cite{Protter}, and for the second-order telegraph equation \cite{Maw1}---\cite{Ortega} are quite different from those for elliptic and parabolic cases.
In the monograph of Protter and  Weinberger \cite{Protter} there was shown that solutions of hyperbolic equations and inequalities do not exhibit the classical formulation of maximum principle. Even in the simplest case of the wave equation in two independent variables $u_{tt}-u_{xx}=0$  the maximum of a nonconstant solution $u = \sin x\sin t$ in a rectangle domain $\{(x,t) :\, x \in [0, \pi],\, t \in [0, \pi] \}$  occurs at an interior point $ \left(\frac{\pi}{2},\,\frac{\pi}{2}\right)$. In Chapter 4 \cite{Protter} the maximum principle for linear second hyperbolic equations of general type, with variable coefficients has also been obtained for Cauchy problems and boundary value problems on characteristics (Goursat problem).

Following R.Ortega, A.Robles-Perez \cite{Ortega}, we introduce the definition of a "weak form" of the maximum principle, which is used for the hyperbolic equations, which will be used later.
\vskip5mm
 {\it Definition 1.} \cite{Ortega}
 Let $ L= Lu$ be linear differential operator, acting on functions $u:\,D\to \mathbb{R}$, in some domain $D$. These functions will belong to the certain family $ B,$ which includes boundary conditions or others requirements. It is said that $ L$ satisfies the maximum principle, if
$$ L\geq 0,\,\,u\in B,$$
implies $u\geq 0$ in $D$.
 \vskip5mm

 In further works of these authors (see \cite{Maw1}, \cite{Maw2},\cite{Maw3}) there was studied  the maximum principle for weak bounded twice periodical solutions from the space $L^{\infty}$ of the telegraph equation with parameter $\lambda$ in lower term, one-, two-, and -tree dimensional spaces, and which includes the cases of variables coefficients. The precise condition for $\lambda$ under which the maximum principle still valid was font.  There was also introduced a method of upper and lower solutions  associated with the nonlinear equation, which allows to obtain the analogous results (uniqueness, existence and regularity theorems) for the  telegraph equations with external nonlinear forcing, applying maximum principle. There was considered also the case when the external forcing belongs to a certain space of measures.

 The maximum principle for general quasilinear hyperbolic systems with dissipation was proved by Kong De Xing \cite{Kong}. There were given two estimates for the solution to the general quasilinear hyperbolic system and introduced the concept of dissipation (strong dissipation and weak dissipation), then state some maximum principles of quasilinear hyperbolic systems with dissipation. Using the maximum principle there were reproved the existence and uniqueness theorems of the global smooth solution to the Cauchy problem for considered quasilinear hyperbolic system.

 So, the problem to prove the maximum principle for the weak solutions stills more complicated and at that time becomes more interesting in the case of fourth-order hyperbolic equations, especially, in the case of non-classical boundary value problems with weak-regularity data.  There are no results on maximum principle even for model case of linear 2-dimensional fourth-order hyperbolic equations with the constant coefficients and without lower terms. Moreover, we can not use the term of usual traces in the cases of initial data of weak regularity, and we come to the notions of $L-$traces, the traces, which associated with differential operator. Let us remind (see, for example, \cite {Bursk1}), that $L-$ traces exist for the weak solutions from space $L^2$ even in the situations when classical notions of traces does not work for such solutions.

 Interesting interpretation of violation of Fredholm property arises as periodicity of characteristic billiard or John's mapping. In the simple case of second order hyperbolic equations there was proved \cite{John}, \cite{Bursk2} that periodicity of John's algorithm is sufficient for violation of the Fredholm property of the  Dirichlet problem. The analogous result is true for the case of fourth-order hyperbolic equations and will be proved in the present paper.

 Therefore, to establish the maximum principle and to obtain results on uniqueness, existence and regularity, kernel dimensionality, Fredholm property for weak solutions for fourth order hyperbolic operators and boundary value problems for them are very important for the reason of their applications and further study of these problems, and it is the mail goal of the paper.

\section{Statement of the problem and auxiliary definitions }

Let us start to establish the maximum principle for the weak solutions to the Cauchy problem
for the Eq.(\ref{eq:01}) in some admissible planar domain. It is  expected, that in hyperbolic case the  characteristics of the equations play a crucial role.

Let $C_j,\,j=1,\,2,\,3,\,4$ be characteristics,  $\Gamma_0:=\{x_1\in [a,\,b],\,x_2=0\},$ and  define domain $\Omega$ as a domain which is restricted by the characteristics $C_j,\,j=1,\,2,\,3,\,4$ and $\Gamma_0.$ Consider also the following Cauchy problem for the Eq. (\ref{eq:01}) on $\Gamma_0:$
\begin{equation}\label{eq:02}
u|_{\Gamma_0}=\varphi(x),\,u^{\prime}_{\nu}|_{\Gamma_0}=\psi(x),\,u^{\prime\prime}_{\nu\nu}|_{\Gamma_0}=\sigma(x),\,u^{\prime\prime\prime}_{\nu\nu\nu}|_{\Gamma_0}=
\chi(x),
\end{equation}
where $\varphi,\,\psi,\,\sigma$ and $\chi$ are given weak regular functions on $\Gamma_0,$ in general case $ \varphi,\,\psi,\,\sigma,\,\chi\in L^2(\Gamma_0),\,\nu-$ is outer normal of $\Gamma_0.$

{\it Definition 2.} We call a domain $D:=\{(x_1,\,x_2):\,x_1\in (-\infty,\,+\infty),\,x_2>0\}$ in the half-plane $x_2 > 0$ an admissible domain if
it has the property that for each point $C\in D$ the corresponding characteristics domain $\Omega$ is also in $D$.  More generally, $D$ is admissible, if it is the finite or countable union of characteristics 5-angles (in the case of fourth-order equations with constant coefficients, that is existence of 4 different and real characteristics lines).

Establishment of the maximum principle in this situation allows us to obtain a local properties of the solution to Cauchy problem (\ref{eq:01})---(\ref{eq:02}) on the arbitrary interior point $C\in D.$

We will consider the weak solution to the problem (\ref{eq:01})---(\ref{eq:02}) from the $D(L)$, domain of definition of maximal operator, associated with the differential operation $L$ in Eq.(\ref{eq:01}). Following \cite{Bur3},\,\cite{Grubb} and \cite{H1}, we remind the corresponding definitions.

In the bounded domain $\Omega$ we consider the linear differential operation ${\mathcal L}$ of the order $m,\,m\geq 2,$ and formally adjoint ${\mathcal L}^+$:
\begin{equation}\label{eq:03}
{\mathcal L}(D_x)=\sum\limits_{|\alpha|\leq m}a_{\alpha}D^{\alpha},\,{\mathcal L}^+(D_x)=\sum\limits_{|\alpha|\leq m}D^{\alpha}(a_{\alpha}),
\end{equation}
where $\alpha=(\alpha_1,\,\alpha_2,...\alpha_n),\,|\alpha|=\alpha_1+\alpha_2+...+\alpha_n$ is multi-index. Note, that for Eq. (\ref{eq:01}) $n=2,\,m=4.$

{\it Definition 3. Minimum operator.} \cite{Bur3}.
Let us consider the differential operation (\ref{eq:03}) on functions from the space $C_0^{\infty}(\Omega).$ The minimum operator $L_0$ is called the extension of the operation from $C_0^{\infty}(\Omega)$ to the set $D(L_0):=\overline{C_0^{\infty}(\Omega)}.$ The closure is realized in the norm of the graph of operator $L$: $||u||_L^2:=||u||_{L_2(\Omega)}^2+||Lu||^2_{L_2(\Omega)}.$

{\it Definition 4. Maximum operator.} \cite{Bur3}.
The maximum operator $L$ is define as the restriction of the differential operation $  {\mathcal L}(D_x)$ to the set $D(L):=\{u\in L^2(\Omega):\,Lu\in L^2(\Omega)\}.$

{\it Definition 5.} \cite{Bur3}.
The operator $\tilde L$ is define as the extension of the minimum operator $L_0,$ to the set $D(\tilde L):=\overline{C^{\infty}(\bar\Omega)}.$

{\it Definition 6. Regular operator.} \cite{Bur3}.
The maximum operator is called regular if  $D(L)=D(\tilde L).$

It is easy to see, that in the case of operation of the fourth-order (\ref{eq:01}), maximal operator is regular and $D(L)=D(\tilde L)=H^4(\Omega),$ $D(L_0)=\stackrel{0}{H^4}(\Omega),$ the  Hilbert Sobolev space of fourthly weak differentiable functions from $L^2(\Omega)$.

The definition of a weak solution to the problem (\ref{eq:01})---(\ref{eq:02}) from the space $D(L)$ is closely connected with the notion of $L-$ traces, that is traces, which are associated with the differential operator $L$.

{\it Definition 7. L-traces.} \cite{Bur5}.
Assume, that for a function $u\in D(\tilde L)$ there exist linear continuous functionals $L_{(p)}u$ over the space $H^{m-p-1/2}(\partial\Omega),\,p=0,1,2...,m-1$, such that the following equality is satisfied:
\begin{equation}\label{eq:04}
(Lu,v)_{L^2(\Omega)}-(u,L^+v)_{L^2(\Omega)}=\sum\limits_{j=0}^{m-1}(L_{(m-1-j)}u,\,\partial^{(j)}_{\nu}v).
\end{equation}

The functionals $L_{(p)}u$ is called the $L_{(p)}-$ traces of the function $u\in D(\tilde L).$ Here $(\cdot,\,\cdot)_{L^2(\Omega)}$ is a scalar product in Hilbert space $L^2(\Omega)$.

For $L^2-$ solutions the notion of $L_{(p)}-$ traces can be realized by the following way.

{\it Definition 8. L-traces.}

Finally, we are going to the definition of the weak solution to the problem (\ref{eq:01})---(\ref{eq:02}):

{\it Definition 9.} We will call the function $u\in D(L)$ a weak solution to the Cauchy problem (\ref{eq:01})--(\ref{eq:02}), if it satisfies to the following integral identity
\begin{equation}\label{eq:05}
(f,\,v)_{L^2(\Omega)}-(u,\,L^+v)_{L^2(\Omega)}=\sum\limits_{j=0}^{3}(L_{(3-j)}u,\,\partial^{(j)}_{\nu}v),
\end{equation}
for any functions $v\in C_0^{\infty}(\Omega).$
The functionals $L_{(p)}u$ is called the $L_{(p)}-$ traces of the function $u$, $p=0,\,1,\,2,\,3,$ and completely determined by the initial functions $\varphi,\,\psi,\,\sigma,\,\chi$ by the following way:
\begin{equation*}
L_{(0)}u=-L(x)u|_{\partial\Omega}=-L(\nu)\varphi;
\end{equation*}
\begin{equation*}
L_{(1)}u=L(\nu)\psi+\alpha_1\varphi_{\tau}^{\prime}+\alpha_2\varphi;
\end{equation*}
\begin{equation}\label{eq:06}
L_{(2)}u=-L(\nu)\sigma+\beta_1\psi_{\tau}^{\prime}+\beta_2\psi+\beta_3\varphi_{\tau\tau}^{\prime\prime}+\beta_4\varphi_{\tau}^{\prime}+\beta_5\varphi;
\end{equation}
\begin{equation*}
L_{(3)}u=L(\nu)\chi+\delta_1\varphi_{\tau\tau\tau}^{\prime\prime\prime}+\delta_2\sigma+\delta_3\psi_{\tau\tau}^{\prime\prime}+\delta_4\psi_{\tau}^{\prime}+
\delta_5\psi+\delta_6\varphi_{\tau\tau}^{\prime\prime}+\delta_7\varphi_{\tau}^{\prime}+\delta_8\varphi.
\end{equation*}
Here $\alpha_i,\,i=1,\,2,\,\beta_j,\,j=1,\,2,...,\,5,$ and $\delta_k,\,k=1,\,...,\,9$ are smooth functions, completely determined by the coefficients of the Eq.(\ref{eq:01}).

{\it Remark 1.} We can use a general form of the operators $\gamma_j$ in the left-hand side of the identity (\ref{eq:05}) instead of operators of differentiation $\partial^{(j)}_{\nu}v$. Indeed, we define $\gamma_j=p_j\gamma,$ where $\gamma: \,u\in H^m(\Omega)\to (u|_{\partial\Omega},\,...,\,u_{\nu}^{(m-1)}|_{\partial\Omega})\in H^{(m)}=H^{m-1/2}(\partial\Omega)\times H^{m-3/2}(\partial\Omega)\times...\times H^{1/2}(\partial\Omega),$ and $p_j:\,H^{(m)}\to H^{m-j-1/2}(\partial\Omega)-$ projection.

As it has been mentioned above,  some examples show (see \cite{Bursk1}) that in the general case the solutions $u\in D(L)$ do
not exist ordinary traces in the sense of distributions even for the simplest hyperbolic equations.
 Indeed, for the wave equation $Lu =\frac{\partial^2 u}{\partial x_1\partial x_2} = 0$ in the unit disk $K:\, |x|=1$, the solution $u(x) = (1-x_1^2)^{-\frac{5}{2}}$
belongs to $L^2(K)$, but $<u|_{\partial K}, 1>_{\partial K} = \infty$ it means that $\lim_{r\to 1-0}\int\limits_{|x|=r}u(x)ds_x=\infty$, such a way
the trace $u|_{\partial K}$ does not exist even as a distribution. However, for every solution $u\in L^2(K)$  the $L_{(0)}-$trace $L_{(0)}u := -L(x) u(x)|_{|x|=1}=-x_1x_2 u(x)|_{|x|=1}\in L^2(\partial K).$  Likewise, $L_{(1)}-$ trace $L_{(1)}u$ exists for every $u\in L^2(K)$:
$$
L_{(1)}u=\left(L(x)u^{\prime}_{\nu}+L^{\prime}_{\tau}u^{\prime}_{\tau}+\frac{1}{2}L^{\prime\prime}_{\tau\tau}u\right)|_{\partial K}\in H^{-\frac{3}{2}}(\partial K).
$$
where $\tau$ is the angular coordinate and $u^{\prime}_{\tau}$
is the tangential derivative, and $L(x)=x_1x_2-$ symbol of the operator $L=\frac{\partial^2 }{\partial x_1\partial x_2}$.

\section{Maximum principle for the weak solutions of Cauchy problem. Existence, uniqueness and regularity of solution}
We prove here the first simple case: the maximum principle for the weak solution of the Cauchy problem (\ref{eq:01})---(\ref{eq:02}) in admissible plane domain $\Omega,$ restricted by the different and  non-congruent characteristics $C_j,\,j=1,\,2,...,\,4$ and initial line $\Gamma_0$.

{\it Theorem 1. Maximum principle.} Let $u\in D(L)$ satisfy the following inequalities:

\begin{equation}\label{eq:07}
 Lu=f\leq 0,\,\,\,x\in D,
 \end{equation}
and
 \begin{equation}
L_{(0)}u\mid_{\Gamma_0}\geq 0,\,L_{(1)}u|_{\Gamma_0}\geq 0,\,L_{(2)}u|_{\Gamma_0}\geq 0,\,
L_{(3)}u|_{\Gamma_0}\geq 0,\label{eq:08}
\end{equation}
then
$u\leq 0$ in $D.$
\vskip3.5mm
{\it Proof.}

Due to homogeneity of the symbol in Eq. (\ref{eq:01}), $L(\xi)=a_0\xi_1^4+a_1\xi_1^3\xi_2+a_2\xi_1^2\xi_2^2+a_3\xi_1\xi_2^3+a_4\xi_2^4=$ \\ $<\xi,\,a^1><\xi,\,a^2><\xi,\,a^3><\xi,\,a^4>,\,\xi=(\xi_1,\,\xi_2)\in \mathbb R^2,$ we can rewrite this equation in the following form:
\begin{equation}
<\nabla,\,a^1><\nabla,\,a^2><\nabla,\,a^3><\nabla,\,a^4>u=f(x).\label{eq:09}
\end{equation}
The vectors $a^j=(a_1^j,\,a_2^j),\,j=1,\,2,\,3,\,4$ are determined by the coefficients $a_i,\,i=0,\,1,\,2,\,3,\,4, $ and $<a,\,b>=a_1\bar b_1+a_2\bar b_2$ is a scalar product in $\mathbb C^2$. It is easy to see that vector $a^j$ is a tangent vector of $j-$th characteristic, slope $\varphi_j$ of which is determined by $-\tan\varphi_j=\lambda_j,\,j=1,\,2,\,3,\,4.$ In what follows, we also consider the vectors $\tilde a^j=(-\bar a_2^j,\,\bar a_1^j),\,j=1,\,2,\,3,\,4.$ It is obvious that $<\tilde a^j,\,a^j>=0,$ so $\tilde a^j$ is a normal vector of $j-$th characteristic.

Use the definitions 7 and 9 for the case  $m=4,$ that is fourth-order operator in Eq. (\ref{eq:01}), and domain  $\Omega, $ which is restricted by the characteristics $C_j,\,j=1,\,2,\,3,\,4$ and $\Gamma_0:$
$$
\int\limits_{\Omega}\{Lu\cdot \bar v-u\cdot\overline{L^+v}\}dx=\sum\limits_{k=0}^{3}\int\limits_{\partial\Omega}L_{(3-k)}u\cdot \partial^{(k)}_{\nu}v\,ds=
$$
$$
=\sum\limits_{k=0}^{3}\int\limits_{C_1}L_{(3-k)}u\cdot \partial^{(k)}_{\nu}v\,ds+\sum\limits_{k=0}^{3}\int\limits_{C_2}L_{(3-k)}u\cdot \partial^{(k)}_{\nu}v\,ds+\sum\limits_{k=0}^{3}\int\limits_{C_3}L_{(3-k)}u\cdot \partial^{(k)}_{\nu}v\,ds+\sum\limits_{k=0}^{3}\int\limits_{C_4}L_{(3-k)}u\cdot \partial^{(k)}_{\nu}v\,ds+
$$
\begin{equation}
+\sum\limits_{k=0}^{3}\int\limits_{\Gamma_0}L_{(3-k)}u\cdot \partial^{(k)}_{\nu}v\,ds.\label{eq:10}
\end{equation}
Using the representation (\ref{eq:09}), we arrive to
$$\int\limits_{\Omega}Lu\cdot \bar v\,dx=\int\limits_{\Omega}<\nabla,\,a^1><\nabla,\,a^2><\nabla,\,a^3><\nabla,\,a^4>u\cdot\bar v\,dx=$$
$$\int\limits_{\partial\Omega}<\nu,\,a^1>\cdot<\nabla,\,a^2><\nabla,\,a^3><\nabla,\,a^4>u\cdot\bar v\,ds-$$
$$\int\limits_{\Omega}<\nabla,\,a^2><\nabla,\,a^3><\nabla,\,a^4>u\cdot\overline{<\nabla,\,a^1>v}\,dx.$$
Integrating by parts further, we obtain:
$$\int\limits_{\Omega}Lu\cdot \bar v\,dx=
\int\limits_{\partial\Omega}<\nu,\,a^1>\cdot<\nabla,\,a^2><\nabla,\,a^3><\nabla,\,a^4>u\cdot\bar v\,ds-$$
$$\int\limits_{\partial\Omega}<\nu,\,a^2>\cdot<\nabla,\,a^3><\nabla,\,a^4>u\cdot\overline{<\nabla,\,a^1>v}\,ds+$$
$$+\int\limits_{\partial\Omega}<\nu,\,a^3>\cdot<\nabla,\,a^4>u\cdot\overline{<\nabla,\,a^2><\nabla,\,a^1>v}\,ds-$$
$$-\int\limits_{\partial\Omega}<\nu,\,a^4>\cdot u\cdot\overline{<\nabla,\,a^3><\nabla,\,a^2><\nabla,\,a^1>v}\,ds+$$
$$+\int\limits_{\Omega}u\cdot\overline{<\nabla,\,a^4><\nabla,\,a^3><\nabla,\,a^2><\nabla,\,a^1>v}\,dx.$$
Since $<\nabla,\,a^4><\nabla,\,a^3><\nabla,\,a^2><\nabla,\,a^1>v=L^+v,$ and determining
$$\tilde L_{(0)}u:=<\nu,\,a^4>u,\,\,\tilde L_{(1)}u:=<\nu,\,a^3>\cdot<\nabla,\,a^4>u,$$
$$\tilde L_{(2)}u:=<\nu,\,a^2>\cdot<\nabla,\,a^3><\nabla,\,a^4>u,$$
 $$\tilde L_{(3)}u=L_{(3)}u=<\nu,\,a^1>\cdot<\nabla,\,a^2><\nabla,\,a^3><\nabla,\,a^4>u$$ that are
analogues of $L-$traces from the formula (\ref{eq:10}). Such a way we have
$$\int\limits_{\Omega}\{Lu\cdot \bar v-u\cdot\overline{L^+v}\}\,dx=
\int\limits_{\partial\Omega}L_{(3)}u\cdot\bar v\,ds-\int\limits_{\partial\Omega}\tilde L_{(2)}u\cdot\overline{<\nabla,\,a^1>v}\,ds+$$
\begin{equation}\label{eq:11}
+\int\limits_{\partial\Omega}\tilde L_{(1)}u\cdot\overline{<\nabla,\,a^2><\nabla,\,a^1>v}\,ds
-\int\limits_{\partial\Omega}\tilde L_{(0)} u\cdot\overline{<\nabla,\,a^3><\nabla,\,a^2><\nabla,\,a^1>v}\,ds.
\end{equation}
Difference between formulas (\ref{eq:10}) and (\ref{eq:11}) is that natural $L_{(3-k)}$ traces in (\ref{eq:10}) are multiplied by $k-$ derivative by outer normal $\nu$ of truncated function $v:\,\partial^{(k)}_{\nu}v,$ on the other hand, in (\ref{eq:11}) we determined by $\tilde L_{(3-k)}$ some expressions which multiplied by differential operators $L_{k}^+v$ of order $k$ and which can serve as  analogous of natural  $L_{(3-k)}$ traces, $k=0,\,1,\,2,\,3.$ So, in the (\ref{eq:11})
$$L^+_1v:=<\nabla,\,a^1>v,\,L^+_2v:=<\nabla,\,a^2><\nabla,\,a^1>v,$$
$$L_0^+v=v,\,L^+_3v:=<\nabla,\,a^3><\nabla,\,a^2><\nabla,\,a^1>v.$$
Let  $v\in Ker L^+$ in (\ref{eq:11}) and calculate $L-$ traces on $\partial\Omega=C_1\cup C_2\cup C_3\cup C_4\cup\Gamma_0$. For instance, for $L_{(3)}u$ we obtain: $L_{(3)}u=<\nu,\,a^1><\nabla,\,a^2><\nabla,\,a^3><\nabla,\,a^4>u, $ and use that $<\nabla,\,a^j>u=$\\$<\nu,\,a^j>u^{\prime}_{\nu}+<\tau,\,a^j>u^{\prime}_{\tau},\,j=1,\,2,\,3,\,4,$ where $\nu-$ normal vector, $\tau-$ tangent vector. Due to presence the product $<\nu,\,a^1>,$  $L_{(3)}u=0$ on characteristic $C_1,$ normal vector $\tilde a^1$ of which is orthogonal to the vector $a^1$. On the other parts of $\partial\Omega$ there will vanish the terms containing $<\nu,a^j>$ on $C_j$.   After that\
$$\int\limits_{\partial\Omega}<\nu,\,a^1><\nabla,\,a^2><\nabla,\,a^3><\nabla,\,a^4>u=\int\limits_{\Gamma_0}L_{(3)}u\,ds+$$
$$<\tilde a^2,\,a^1><a^2,\,a^2><\tilde a^2,\,a^3><\tilde a^2,\,a^4>\int\limits_{C_2}u_{\nu\nu\tau}\,ds+$$
$$<\tilde a^3,\,a^1><\tilde a^3,\,a^2>< a^3,\,a^3><\tilde a^3,\,a^4>\int\limits_{C_3}u_{\nu\nu\tau}\,ds+$$
$$<\tilde a^4,\,a^1><\tilde a^4,\,a^2><\tilde a^4,\,a^3>< a^4,\,a^4>\int\limits_{C_4}u_{\nu\nu\tau}\,ds+$$
$$\left\{<\tilde a^2,\,a^1><a^2,\,a^2><\tilde a^2,\,a^3><a^2,\, a^4>+<\tilde a^2,\,a^1><a^2,\,a^2>< a^2,\,a^3><\tilde a^2,\, a^4>\right\}\int\limits_{C_2}u_{\tau\tau\nu}\,ds+$$
$$ \left\{<\tilde a^3,\,a^1><\tilde a^3,\,a^2><a^3,\,a^3><a^3,\,a^4>+<\tilde a^3,\,a^1><a^3,\,a^2>< a^3,\,a^3><\tilde a^3,\, a^4>\right\}\int\limits_{C_3}u_{\tau\tau\nu}\,ds+$$
$$\left\{<\tilde a^4,\,a^1><\tilde a^4,\,a^2><a^4,\,a^3><a^4,\, a^4>+<\tilde a^4,\,a^1><a^4,\,a^2>< \tilde a^4,\,a^3><a^4,\, a^4>\right\}\int\limits_{C_4}u_{\tau\tau\nu}\,ds+$$
$$<\tilde a^2,\,a^1><a^2,\,a^2>< a^2,\,a^3><a^2,\,a^4>\int\limits_{C_2}u_{\tau\tau\tau}\,ds+$$
$$<\tilde a^3,\,a^1><a^3,\,a^2>< a^3,\,a^3><a^3,\,a^4>\int\limits_{C_3}u_{\tau\tau\tau}\,ds+$$
$$<\tilde a^4,\,a^1><a^4,\,a^2>< a^4,\,a^3><a^4,\,a^4>\int\limits_{C_4}u_{\tau\tau\tau}\,ds+
\alpha_{4,1}\int\limits_{C_2}u_{\nu\nu}\,ds+\alpha_{4,2}\int\limits_{C_3}u_{\nu\nu}\,ds+\alpha_{4,3}\int\limits_{C_4}u_{\nu\nu}\,ds+$$
$$\alpha_{5,1}\int\limits_{C_2}u_{\nu\tau}\,ds+\alpha_{5,2}\int\limits_{C_3}u_{\nu\tau}\,ds+\alpha_{5,3}\int\limits_{C_4}u_{\nu\tau}\,ds
+\alpha_{6,1}\int\limits_{C_2}u_{\tau\tau}\,ds+\alpha_{6,2}\int\limits_{C_3}u_{\tau\tau}\,ds+\alpha_{6,3}\int\limits_{C_4}u_{\tau\tau}\,ds+$$
$$\alpha_{7,1}\int\limits_{C_2}u_{\nu}\,ds+\alpha_{7,2}\int\limits_{C_3}u_{\nu}\,ds+\alpha_{7,3}\int\limits_{C_4}u_{\nu}\,ds+
\alpha_{8,1}\int\limits_{C_2}u_{\tau}\,ds+\alpha_{8,2}\int\limits_{C_3}u_{\tau}\,ds+\alpha_{8,3}\int\limits_{C_4}u_{\tau}\,ds.$$
Here correspondent coefficients $\alpha_{i,j}$ were numerated as follows: first index $i$ indicates the derivative of $u$:  $1)\,u_{\nu\nu\tau},\,2)\,u_{\nu\tau\tau},\,3)\,u_{\tau\tau\tau},\,4)\, u_{\nu\nu},\,5)\,u_{\nu\tau},\,6)\,u_{\tau\tau},\,7)\, u_{\nu}\, 8)\,u_{\tau},$ the second index $j$ indicates $j+1-$th characteristic, $j=1,\,2,\,3.$
So, now the  formula (\ref{eq:10}) has the form:
$$\int\limits_{\Omega}Lu\,dx=\int\limits_{\Gamma_0}L_{(3)}u\,ds+
\alpha_{1,1}\int\limits_{C_2}u_{\nu\nu\tau}\,ds+\alpha_{1,2}\int\limits_{C_3}u_{\nu\nu\tau}\,ds+\alpha_{1,3}\int\limits_{C_4}u_{\nu\nu\tau}\,ds+$$
$$\alpha_{2,1}\int\limits_{C_2}u_{\tau\tau\nu}\,ds+\alpha_{2,2}\int\limits_{C_3}u_{\tau\tau\nu}\,ds+\alpha_{2,3}\int\limits_{C_4}u_{\tau\tau\nu}\,ds+$$
$$\alpha_{3,1}\int\limits_{C_2}u_{\tau\tau\tau}\,ds+\alpha_{3,2}\int\limits_{C_3}u_{\tau\tau\tau}\,ds+\alpha_{3,3}\int\limits_{C_4}u_{\tau\tau\tau}\,ds+$$
$$\alpha_{4,1}\int\limits_{C_2}u_{\nu\nu}\,ds+\alpha_{4,2}\int\limits_{C_3}u_{\nu\nu}\,ds+\alpha_{4,3}\int\limits_{C_4}u_{\nu\nu}\,ds+$$
$$\alpha_{5,1}\int\limits_{C_2}u_{\nu\tau}\,ds+\alpha_{5,2}\int\limits_{C_3}u_{\nu\tau}\,ds+\alpha_{5,3}\int\limits_{C_4}u_{\nu\tau}\,ds+
\alpha_{6,1}\int\limits_{C_2}u_{\tau\tau}\,ds+\alpha_{6,2}\int\limits_{C_3}u_{\tau\tau}\,ds+\alpha_{6,3}\int\limits_{C_4}u_{\tau\tau}\,ds+$$
$$\alpha_{7,1}\int\limits_{C_2}u_{\nu}\,ds+\alpha_{7,2}\int\limits_{C_3}u_{\nu}\,ds+\alpha_{7,3}\int\limits_{C_4}u_{\nu}\,ds+
\alpha_{8,1}\int\limits_{C_2}u_{\tau}\,ds+\alpha_{8,2}\int\limits_{C_3}u_{\tau}\,ds+\alpha_{8,3}\int\limits_{C_4}u_{\tau}\,ds.$$
Coefficients $\alpha_{i,j}$ are constant and depend on only from Eq. (\ref{eq:01}) coefficients $a_0,\,a_1,\,a_2,\,a_3,\,a_4.$ By analogous way we calculate others $L-$ traces, $L_{(0)}u,\,L_{(1)}u$ and $L_{(2)}u.$

To obtain the statement of the Theorem 1, we choose some arbitrary point $C\in D$ in admissible plane domain $D,$  draw through this point two arbitrary characteristics, $C_1$ and $C_2$. Another two characteristics ($C_3$ and $C_4$) we draw through the ends $a$ and $b$ of initial line $\Gamma_0$. We determine the points $O_1$ and $O_2$ as intersections of $C_1,\,C_3$ and $C_2,\,C_4$ correspondingly: $O_1=C_1\cap C_3,\,O_2=C_2\cap C_4.$ Such a way, domain $\Omega$ is a pentagon $aO_1CO_2b.$ The value of the function $u$ at the point $C\in D,\, u(C)$ we estimate from the last equality, integrating by the characteristics $C_1$ and $C_2$ and using conditions (\ref{eq:02}), (\ref{eq:06})---(\ref{eq:08}). Since, the chosen point $C\in D$ is arbitrary, we arrive at $u\leq 0$ in $D.$

{\it Remark 2.} In the case of classical solution of the Cauchy problem for the second order hyperbolic equations of the general form with the constant coefficients the statement of the Theorem 1 coincides with the result of \cite{Protter}.  In this case conditions (\ref{eq:08}) have usual form without using the notion of $L-$traces (see \cite{Protter}):
\begin{equation*}
u|_{\Gamma_0}\leq 0,\,\,u^{\prime}_{\nu}|_{\Gamma_0}\leq 0.
\end{equation*}

\section{Method of  equation-domain duality and its application to the Goursat problem}

We consider here the method of equation–domain duality (see also in [\cite{Bur3}],\,[\cite{Bursk1}] ) for the study of the Goursat  problem. This method allows us to reduce the Cauchy problem (\ref{eq:01})---(\ref{eq:02}) in bounded domain $\Omega$  to equivalent the Goursat boundary value problem. We will show that the method of equation–domain duality can be applied also to boundary value problems in the generalized statement, and we will extend the previously obtained results of the study the Goursat problem for specific equations with constant coefficients.
First of all we consider the method of equation-domain duality for the case of classical, smooth solutions.
\subsection{The method of equation-domain duality for the case of classical, smooth solutions.}
Let $\Omega\in \mathbb R^n$ be a bounded domain defined by the inequality $P(x) > 0$ with some real  polynomial $P(x)$. The equation $P(x)=0$ denotes the boundary $\partial\Omega$. It is assumed that the boundary of the domain is non degenerate for $P$, i.e., $|\nabla P|\neq= 0$ on $\partial\Omega$. Consider the general boundary value problem with $\gamma$ conditions on the boundary for the differential operator $L$ (\ref{eq:03}) of the order $m$, and $\gamma\leq m$:
\begin{equation}\label{eq:12}
L(D_x)u=f(x),\,u|_{\partial\Omega}=0,\,u^{\prime}_{\nu}|_{\partial\Omega}=0,\,...,\,u_{\nu}^{(\gamma-1)}|_{\partial\Omega}=0.
\end{equation}
 By the equation-domain duality we mean (see [\cite{Bur3}]) a correspondence (in the sense of Fourier transform)  between problem (\ref{eq:12}) and the equation
 \begin{equation}\label{eq:13}
 P^{m-\gamma}(-D_{\xi})\{L (\xi)w(\xi)\} = \hat f(\xi).
  \end{equation}
  This correspondence is described in the following lemma.

  {\bf Lemma 1.} For any nontrivial solution of problem (\ref{eq:12}) in the space of smooth functions $C^m(\bar\Omega)$, there exists a nontrivial analytic solution $w$ of equation (\ref{eq:13})  from the space $\mathbb C^n$  in a class $Z_{\Omega}^m$ of entire functions. The class $Z_{\Omega}^m$ is defined as the space of Fourier transforms of functions of the form $\theta_{\Omega}\eta$, where $\eta\in C^m(\mathbb R^n)$, and $\theta_{\Omega}$ is the characteristic function of the domain $\Omega$, $w(\xi) = \widehat {\theta_{\Omega}u}.$ The function $f(x)$ is assumed to be extended by zero beyond the boundary of the domain.

  {\it Proof.} Let $m=4,\,\gamma=2$ and consider the following Dirichlet problem for the  fourth-order operator in (\ref{eq:01}):
 \begin{equation}\label{eq:12.5}
L(D_x)u=f,\,u|_{P(x)=0}=f,\,u^{\prime}_{\nu}|_{P(x)=0}=0.
\end{equation}
  Let also $u\in C^4(\bar\Omega)$ be a classical solution to the problem (\ref{eq:12.5}. Denote by $\tilde u\in C^4(\mathbb R^2)$ the extension of $u,$  and apply the fourth-order operator $L(D_x)$ in (\ref{eq:01}) to the product $\tilde u\theta_{\Omega},$ where $\theta_{\Omega}$ is a characteristics function of the domain $\Omega$: $\theta_{\Omega}=1$ in $\Omega$ and $\theta_{\Omega}=0$ out of $\Omega.$
  We obtain:
  $$ L(D_x)(\tilde u\theta_{\Omega})=\theta_{\Omega}L(D_x)\tilde u+\tilde u L(D_x)\theta_{\Omega}+$$
  $$+L^{(1)}_3(D_x)\tilde u<\nabla,\,a^1>\theta_{\Omega}+L^{(2)}_3(D_x)\tilde u<\nabla,\,a^2>\theta_{\Omega}+$$
  $$+L^{(3)}_3(D_x)\tilde u<\nabla,\,a^3>\theta_{\Omega}+L^{(4)}_3(D_x)\tilde u<\nabla,\,a^4>\theta_{\Omega}+$$
  $$+L^{(1)}_3(D_x)\theta_{\Omega}<\nabla,\,a^1>\tilde u+L^{(2)}_3(D_x)\theta_{\Omega}<\nabla,\,a^2>\tilde u+$$
  $$+L^{(3)}_3(D_x)\theta_{\Omega}<\nabla,\,a^3>\tilde u+L^{(4)}_3(D_x)\theta_{\Omega}<\nabla,\,a^4>\tilde u+$$
  $$+L^{(1,2)}_2(D_x)\tilde u<\nabla,\,a^1><\nabla,\,a^2>\theta_{\Omega}+L^{(1,3)}_2(D_x)\tilde u<\nabla,\,a^1><\nabla,\,a^3>\theta_{\Omega}+$$
  $$+L^{(1,4)}_2(D_x)\tilde u <\nabla,\,a^1><\nabla,\,a^4>\theta_{\Omega}+
   L^{(2,3)}_2(D_x)\tilde u<\nabla,\,a^2><\nabla,\,a^3>\theta_{\Omega}+$$
   $$+L^{(2,4)}_2(D_x)\tilde u<\nabla,\,a^2><\nabla,\,a^4>\theta_{\Omega}+L^{(3,4)}_2(D_x)\tilde u
  <\nabla,\,a^3><\nabla,\,a^4>\theta_{\Omega}+$$
  $$+L^{(1,2)}_2(D_x)\theta_{\Omega}<\nabla,\,a^1><\nabla,\,a^2>\tilde u+L^{(1,3)}_2(D_x)\theta_{\Omega}<\nabla,\,a^1><\nabla,\,a^3>\tilde u+$$ $$+L^{(1,4)}_2(D_x)\theta_{\Omega} <\nabla,\,a^1><\nabla,\,a^4>\tilde u+
  L^{(2,3)}_2(D_x)\theta_{\Omega}<\nabla,\,a^2><\nabla,\,a^3>\tilde u+$$
  $$+L^{(2,4)}_2(D_x)\theta_{\Omega}<\nabla,\,a^2><\nabla,\,a^4>\tilde u+L^{(3,4)}_2(D_x)\theta_{\Omega}
  <\nabla,\,a^3><\nabla,\,a^4>\tilde u.$$
  Here $L^{(j)}_3(D_x)$ and $L^{(j,k)}_2(D_x),\,j,\,k=1,2,3,4$ are some  differential operations of the $3-$ and $2-$ order correspondingly, defined by the fourth-order differential operation $L(D_x)$ in (\ref{eq:01}):
  $$L^{(j)}_3(D_x)=\frac{L(D_x)}{<\nabla,\,a^j>},\,j=1,...,\,4,$$
  $$L^{(j,k)}_2(D_x)=\frac{L(D_x)}{<\nabla,\,a^j><\nabla,\,a^k>},\,j\neq k,\,j,k=1,..\,4.$$
  Since $\tilde u$ is a solution of the equation (\ref{eq:01}), we arrive to
  \begin{equation}\label{eq:13.5}
  L(D_x)(\tilde u\theta_{\Omega})=\theta_{\Omega}f+\tilde u L(D_x)\theta_{\Omega}+A^{(1)}(x)(\delta_{\partial\Omega})^{\prime\prime}_{\nu\nu}+A^{(2)}(x)(\delta_{\partial\Omega})^{\prime}_{\nu}+
  A^{(3)}(x)\delta_{\partial\Omega},
  \end{equation}
  where $A^{(j)}(x)-$ are some smooth functions, depend on coefficients $a^k,\,k=1,...,\,4$ of the equation (\ref{eq:01}) and $j-$ derivatives of function $u$ by outer normal $\nu:$ $u_{\nu}^{(j)},$ and tangent direction $\tau$: $u_{\tau}^{(j)},\,j=1,2,3.$
  Taking into account the conditions (\ref{eq:12.5}), and $<(\delta_{\partial\Omega})^{\prime}_{\nu},\phi>=-<\delta_{\partial\Omega},\phi^{\prime}_{\nu}>=- \int\limits_{\partial\Omega}\bar{\phi^{\prime}_{\nu}}(s)\,ds,\,\forall\,\psi\in \mathcal D(\mathbb R^2),$ we have that $\tilde u L(D_x)\theta_{\Omega}+A^{(1)}(x)(\delta_{\partial\Omega})^{\prime\prime}_{\nu\nu}=0$, and
  $A^{(2)}(x)(\delta_{\partial\Omega})^{\prime}_{\nu}=-\int\limits_{\partial\Omega}(A^{(2)}(s))^{\prime}_{\nu}\,ds=\tilde A^{(3)}(x)\delta_{\partial\Omega},$
  from (\ref{eq:13.5}) we obtain
  \begin{equation}\label{eq:13.6}
  L(D_x)(\tilde u\theta_{\Omega})=\theta_{\Omega}f+B^{(3)}(x)\delta_{\partial\Omega},
  \end{equation}
  where $B^{(3)}(x)=\tilde A^{(3)}(x)+A^{(3)}(x)$ is some  smooth function, depend on coefficients $a^k,\,k=1,...,\,4$ of the equation (\ref{eq:01}) and $3-$ derivatives of function $u$ by outer normal $\nu:$ $u_{\nu}^{\prime\prime\prime},$ and tangent direction $\tau$: $u_{\tau}^{\prime\prime\prime}.$
  Let us multiply (\ref{eq:13.6}) by $P^2(x)$, so that $P^2(x)B^{(3)}(x)\delta_{\partial\Omega}=0$, due to $P(x)=0$ on $\partial\Omega$,  and after that we apply the Fourier transform:
  \begin{equation*}
  P^2(-D_{\xi})(v(\xi))=\hat f.
  \end{equation*}
  Here
  $v(\xi)=L(\xi)w(\xi),\,w(\xi)=\widehat{\tilde u\theta_{\Omega}},$ Fourier transform of the function $\tilde u\theta_{\Omega}.$ Such a way we arrive to the dual problem (\ref{eq:13}). Function $w(\xi)\in Z_{\Omega}^4,$ the space of entire functions, which is defined as the space of Fourier transforms of functions of the form $\theta_{\Omega}\tilde u$, where $\tilde\in C^4(\mathbb R^n)$, see, for example, \cite{H1}. The Lemma is proved.

As an application of the Lemma 1, let us consider the Dirichlet problem problem for the fourth-order hyperbolic equation (\ref{eq:01})  in the unit disk $K=\{x\in\mathbb R^2:\,|x|< 1\}:$
   \begin{equation}\label{eq:14}
   u|_{|x|=1}=0,\,u^{\prime}_{\nu}|_{|x|=1}=0.
   \end{equation}
  For this case $m=4,\,\gamma=2,\,m-\gamma=2$ and we arrive to the following dual problem:
  \begin{equation}\label{eq:15}
 \Delta^2 v=\hat f(\xi),\,\,v|_{L(\xi)=0}=0,
  \end{equation}
  $v=L(\xi)w(\xi)$. Taking into account the representation (\ref{eq:09}), the condition $w|_{L(\xi)=0}=0$ equivalent to the following four conditions:
  \begin{equation}\label{eq:16}
  w|_{<\xi,\,a^1>=0}=0,\,w|_{<\xi,\,a^2>=0}=0,\,w|_{<\xi,\,a^3>=0}=0,\,w|_{<\xi,\,a^4>=0}=0.
  \end{equation}
 Since $<\xi,\,a^j>=0$ is a characteristics, $j=1,\,2,...,4$ we conclude, that problem (\ref{eq:15}) is a Goursat problem. Method of equation-domain duality allows us to reduce the problem of solvability of a boundary value problem for a typeless equation (particulary, hyperbolic type) to an analogous problem for an equation, possibly of less complicated structure and of lower order (in particular, for elliptic type equation, see (\ref{eq:15})).
\subsection{The method of equation-domain duality for the case of weak solutions and solutions from $D(L).$}
We prove here the analog of the Lemma 1 for the case of solutions $u\in D(L).$
From the Definition 8 and formulae (\ref{eq:06}), for any function $u\in H^m(\Omega), \,m\geq 4,$ $L_{(p)}u-$ traces can be expressed by the following way: $L_{(p)}u=\sum\limits_{k=0}^p\alpha_{p,k}\partial^k_{\nu}u|_{\partial\Omega},\,p=0,1,2,3.$ For $p=0,$ $L_{(0)}u=u|_{\partial\Omega},$ there coincide $L_{(0)}-$ trace and usual trace.
If $u\in D(L),$ then we consider the following boundary value problem
\begin{equation}\label{eq:17}
L(D_x)u=f(x),\,L_{(0)}u=0,\,L_{(1)}u=0,\,...,\,L_{(\gamma-1)}u=0,\,\gamma\leq m.
\end{equation}
instead of (\ref{eq:12}). For example, in the case of the Dirichlet problem (\ref{eq:12.5}) for the fourth-order operator (\ref{eq:01}), and for $u\in D(L)$ we have
\begin{equation}\label{eq:18}
L(D_x)u=f(x),\,L_{(0)}u=0,\,L_{(1)}u=0,\gamma=2<m=4.
\end{equation}
The equation-domain duality principle for the solutions $u\in D(L)$ of the problem (\ref{eq:17}) is assumed as the correspondence (in the sense of Fourier transform)  between the problem (\ref{eq:17}) and the equation (\ref{eq:13}), realized by the following statement, which is analog of the Lemma 1 for the solutions $u\in D(L):$

{\bf Lemma 2.}  For any nontrivial solution of problem \ref{eq:17} in the space  $D(L)$, there exists a nontrivial analytic solution $w$ of equation (\ref{eq:13})  from the space $\mathbb C^n$  in a class $Z_{\Omega}$ of entire functions. The class $Z_{\Omega}$ is defined as the space of Fourier transforms of functions from the set $V=\{v:\, {\rm there}\,{\rm exists}\, {\rm some}\, {\rm function}\, u\in D(L),\,{\rm such}\,{\rm that}:\,v=u\,{\rm in}\,\Omega,\, v=0,\,{\rm out}\, {\rm of}\,\bar{\Omega}\},$ $w(\xi) = \widehat {v}.$ The function $f(x)$ is assumed to be extended by zero beyond the boundary of the domain.

The proof is following from the Definition 9, substituting into the equality (\ref{eq:05}) $v(x)=P^{m-\gamma}(x)e^{i(x,\tilde a^j)}\in ker(L^+),\,j=1,...,4.$ Function $w(\xi) = \widehat {v}\in Z_{\Omega},$ the space of
entire functions (see, for instance, the Paley-Wiener theorem,\cite{H1}.
\section{Connection between Cauchy and Dirichlet problems. Existence and uniqueness of solutions for the hyperbolic equations.}

The main result of this section is the following theorem on existence and solution uniqueness of the Cauchy problem (\ref{eq:01})--(\ref{eq:02}).

{\bf Theorem 2.} Let us assume that there exist four functions $L_3,\,L_2,\,L_1,\,L_0\in L^2(\partial\Omega),$ satisfying conditions
\begin{equation}\label{eq:22}
\int\limits_{\partial\Omega}\{L_3(x)Q(-\tilde a^j\cdot x)+L_2(x)Q^{\prime}(-\tilde a^j\cdot x)+L_1(x)Q^{\prime\prime}(-\tilde a^j\cdot x)+L_0(x)Q^{\prime\prime\prime}(-\tilde a^j\cdot x)\}dS_x=\int\limits_{\Omega}f(x)\overline{Q(-\tilde a^j\cdot x)}dx,
\end{equation}
for any polynomial $Q\in C[z]$ from the kernel $Ker L^+$ of the operator $L^+,$ $Q(-\tilde a^j\cdot x),\,j=1,\,2,\,3,\,4.$

Then,  there exists a unique solution $u\in D(L)$ of the Cauchy problem (\ref{eq:01})--(\ref{eq:02}), whose $L-$ traces are the given functions  $L_3,\,L_2,\,L_1,\,L_0: $ $L_j=L_{(j)-}$ trace, $j=0,\,1,\,2,\,3,$ and connected with the initial data $\varphi,\,\psi,\sigma,\,\chi$ by the relations (\ref{eq:06}).

{\it Proof.} First of all, we prove the existence of the solution to the Cauchy problem (\ref{eq:01})--(\ref{eq:02}) from the space $D(L).$

Let us consider the auxiliary Dirichlet problem for the properly elliptic eight-order operator $\Delta^4$ with the given boundary conditions $\varphi,\,\psi,\sigma,\,\chi:$
\begin{equation}\label{eq:23}
\Delta^4\omega=0,\,\,\omega|_{\partial\Omega}=\varphi,\,\omega_{\nu}|_{\partial\Omega}=\psi,\,\omega_{\nu\nu}|_{\partial\Omega}=\sigma,\,
\omega_{\nu\nu\nu}|_{\partial\Omega}=\chi.
\end{equation}
It is well known that solution of the problem (\ref{eq:23}) exists and belongs to the space $H^m(\Omega),\,m\geq 4.$
We find the solution $u$ of the Cauchy problem in the form
\begin{equation}\label{eq:24}
u=\omega+v,
\end{equation}
where $v$ is a solution of the following problem with null boundary data:
\begin{equation}\label{eq:25}
L(D_x)v=-L(D_x)\omega+f(x),\,v|_{\partial\Omega}=0,\,v_{\nu}|_{\partial\Omega}=0,\,v_{\nu\nu}|_{\partial\Omega}=0,\,
v_{\nu\nu\nu}|_{\partial\Omega}=0.
\end{equation}
Since the $L-$traces of the function $v$ are zero and operator $L$ is regular, we conclude, that $v\in D(L_0)$ and prove the resolvability of the operator equation with the minimum operator $L_0(D_x):$
\begin{equation}\label{eq:26}
L_0(D_x)v=-L\omega+f(x)
\end{equation}
in the space $D(L_0)$.

For resolvability of the operator equation (\ref{eq:26}) with the minimum operator $L_0(D_x)$ it is necessary and sufficiently that the right-hand part  satisfies the following Fredholm condition
\begin{equation}\label{eq:27}
\int\limits_{\Omega}\{-L\omega+f(x)\}\overline{Q(x)}dx=0,
\end{equation}
for any $Q\in Ker\,L^+.$

We use formula (\ref{eq:04}) for the case of function $\omega$ and fourth-order operator, $m=4,$ and taking into account the boundary conditions (\ref{eq:23}), which  mean that the functions $L_0,\,L_1,\,L_2,\,L_3$ are $L-$ traces for the function $\omega,$ conditions (\ref{eq:22}), we arrive to the fulfilment (\ref{eq:27}) for any $Q\in Ker\,L^+,$ and as consequences, we arrive to the resolvability of Eq.(\ref{eq:26}) in $D(L_0).$ Such a way, taking into account the representation (\ref{eq:24}), we arrive to the conclusion on existence solution $u\in D(L).$

Uniqueness of solution follows from the established above maximum principle for the solutions of the Cauchy problem.
Theorem is proved.

\vskip5mm
{\it Remark 3.} For given boundary data $(L_3,\,L_2,\,L_1,\,L_0)\in H^{m-7/2}(\partial\Omega)\times H^{m-5/2}(\partial\Omega)\times H^{m-3/2}(\partial\Omega)\times H^{m-1/2}(\partial\Omega),\,m\geq 4,$ and $f\in H^{m-4}(\Omega),\,m\geq 4,$ in the case of ellipticity equation (\ref{eq:01}), the solution $u\in H^m(\Omega),\,m\geq 4 $ (see [\cite{Bur5}]). But for hyperbolic case it is not true, because the symbol $L(\xi)$ has four real roots, and using Fourier transform and Lemma 2, we arrive at a decreasing of solution regularity as stated in the Theorem 2.
\vskip5mm
{\it Remark 4.} The problem of resolvability the Cauchy problem (\ref{eq:01})--(\ref{eq:02}) is reduced to the integral moment problem (\ref{eq:22}).

\subsection{The Dirichlet problem.}
In some bounded domain $\Omega\in {\mathbb R}^2$ with elliptic boundary $\partial\Omega=\{x:\,P(x)=0\}$ we consider the following Dirichlet problem for the fourth-order hyperbolic equation (\ref{eq:01}):
\begin{equation}\label{eq:28}
L_{(0)}u|_{P(x)=0}=\varphi,\,\,L_{(1)}u_{\nu}|_{P(x)=0}=\psi.
\end{equation}

Connection between the Dirichlet problem (\ref{eq:01}), (\ref{eq:28}) and corresponding Cauchy problem is assumed as follows: let there exists some solution $u^*\in D(L)$ of the Dirichlet problem (\ref{eq:01}), (\ref{eq:28}), then we can construct $L_{(j)}u^*-$traces, that are functions $L_3,\,L_2,\,L_1,\,L_0$ from the Theorem 2, which are satisfied the condition (\ref{eq:22}). From the Theorem 2 it means that the Cauchy problem is solvable in $D(L)$. To prove the solvability of the Dirichlet problem (\ref{eq:01}), (\ref{eq:28}) in $D(L),$  we have to show that there are exist pair of functions $L_2,\,L_3\in L^2(\partial \Omega),$ which is uniquely determined by the $L_{(0)},\,L_{(1)}-$ traces of the Dirichlet problem (\ref{eq:28}). Such a way we arrive to the following inhomogeneous moment problem of determining unknown functions $L_3,\,L_2$ via known left-hand side:
\begin{equation}\label{eq:29}
\int\limits_{\partial\Omega}\{L_3(x)Q(-\tilde a^j\cdot x)+L_2(x)Q^{\prime}(-\tilde a^j\cdot x)\}dS_x=
\end{equation}
\begin{equation*}
=\int\limits_{\Omega}f(x)\overline{Q(-\tilde a^j\cdot x)}dx-\int\limits_{\partial\Omega}\{L_{(1)}(x)Q^{\prime\prime}(-\tilde a^j\cdot x)+L_{(0)}(x)Q^{\prime\prime\prime}(-\tilde a^j\cdot x)\}dS_x
\end{equation*}
for any polynomial $Q\in C[z]$ from the kernel $Ker L^+$ of the operator $L^+,$ $Q(-\tilde a^j\cdot x),\,j=1,\,2,\,3,\,4.$
Thus, solvability of the Dirichlet problem (\ref{eq:28}) in $D(L)$ reduces to the solvability of moment problem (\ref{eq:29}).

{\bf Theorem 3.} For solvability of the Dirichlet problem (\ref{eq:01}), (\ref{eq:28}) in $D(L)$, it is necessary and sufficient that there is exist some solution $(L_3^*(x),\,L_2^*(x))\in L^2(\partial\Omega)\times L^2(\partial\Omega)$ of the moment problem (\ref{eq:29}). Then $L_3^*(x)=L_{(3)}-$ trace, and $L_2^*(x)=L_{(2)}-$ trace.

{\it Remark 5.} In the particular cases of domain $\Omega$ there can be found the explicit formulas for the evaluation of a couple of functions
$(L_3^*(x),\,L_2^*(x))\in L^2(\partial\Omega)\times L^2(\partial\Omega)$ via the known $L_{(0)},\,L_{(1)}-$ traces. For example, the case of unit disk was considered in \cite{Bur5}.

\section{The role of characteristics billiard for the Fredholm property.}

In this section we consider the case when solution uniqueness of the problems considered above is break down, moreover, the Fredholm property does not hold. In the work \cite{Bur3} there was proved the theorem on Fredholm violation of the Dirichlet problem in $C^m(\Omega),\,m\geq 4$ for the case of typeless PDE. Taking into account the Lemma 2, we arrive to the analogous result in $L^2(\Omega):$

{\bf Theorem 4.} The homogeneous Dirichlet problem $(\ref{eq:01})^0$, $(\ref{eq:28})^0$ has a nontrivial solution in $L^2(\Omega)$ if and only if
\begin{equation}\label{eq:30}
\varphi_j-\varphi_k=\frac{\pi p_{jk}}{q},
\end{equation}
with some $p_{jk},\,q\in {\mathbb Z},\,j,k=1,2,3,4.$
Under conditions (\ref{eq:30}) there exists a countable set of linearly independent polynomial solutions in the form:
\begin{equation}\label{eq:31}
u(x)=\sum\limits_{j=1}^4C_j\left(\frac{1}{2q}T_{q}(-\tilde a^j\cdot x)-\frac{1}{2(q-2)}T_{q-2}(-\tilde a^j\cdot x)    \right).
\end{equation}
Here $T_{q}(-\tilde a^j\cdot x)$ are Chebyshev's polynomials, and $\frac{1}{2q}T_{q}(-\tilde a^j\cdot x)-\frac{1}{2(q-2)}T_{q-2}(-\tilde a^j\cdot x)
\in Ker L^+, \,j=1,2,3,4.$

The necessity of condition (\ref{eq:30}) follows from the equation-domain duality (in the case of unit disk), see Lemma 2; sufficiency is proved by the construction of nontrivial polynomial solutions  (\ref{eq:31}). It is remarkable the fact, that the theorem 4 is true for all types of operator $L$. Here we discuss the conditions (\ref{eq:30}) for the hyperbolic case, in which these conditions mean the periodicity of characteristics billiard or John's mapping.

\subsection{Characteristic billiard.}

  For the domain $\Omega,$ which is convex with respect to characteristics, we construct the mappings $T_j,\,j=1,...,4$ for the fourth order hyperbolic equations by the following way. Let $M_j$ be some point on $\partial\Omega.$ Passing through the point $M_j$ $j-$th characteristic, which angle of slope is $\varphi_j,$ we obtain some point $M_{j+1}\in\partial\Omega.$ Such a way, $T_j$ is a mapping which transforms $M_j$ into $M_{j+1}$ on the $j-$characteristic direction with $\varphi_j$ angle of slope, $j=1,2,3,4.$ We apply the mapping $T_1$ for the point $M_1\in\partial\Omega$ and obtain the point $M_2$. After that we apply the mapping $T_2$ for the point $M_2$ and obtain the point $M_3.$ We transform $M_3$ into $M_4$ on direction of characteristic, which angle of slope equals $\varphi_3$, and finally, we transform $M_4$ into $M_5$ on direction of fourth characteristic. Denote by $T=T_4\circ T_3\circ T_2\circ T_1: M_1\in\partial\Omega\to M_5\in\partial\Omega,$ $T$ is called the John's mapping. Characteristic billiard is understood as a discrete dynamical system on $\partial\Omega$, i.e., an action of the group ${\mathbb Z}.$

Some point $M\in\partial\Omega$ is called a periodic point, if there exists some $n\in{\mathbb N}$ such that $T^n(M)=M.$ Minimal $n$ for which condition $T^n(M)=M$ holds is called the period of the point $M.$  In the case of second order hyperbolic equations there was proved [\cite{Bursk2}] that periodicity of John's algorithm is sufficient for violation of the Fredholm property of the  Dirichlet problem. The analogous result is true for the case of fourth-order hyperbolic equations (\ref{eq:01}). Let us consider the model case of domain $\Omega=K,$ unit disk in ${\mathbb R}^2$.

Let us show that conditions (\ref{eq:30}) are necessary and sufficient for the periodicity of John's algorithm.
It is cleat that
\begin{equation}\label{eq:32}
T_j(M(\tau))=2\varphi_j-\tau,
\end{equation}
where $\tau$ is angular parameter of the point $M\in K.$ From (\ref{eq:32}) it follows
$$T^n(M)=2n(\varphi_4-\varphi_3+\varphi_2-\varphi_1)+\tau=2n(\varphi_4-\varphi_3+\varphi_2-\varphi_1)+2\pi m+\tau,$$
for any $m\in{\mathbb Z.}$ Under conditions (\ref{eq:30}) any point $M\in K$ is periodical, thus, the John's algorithm is periodical. If now mapping $T$ is periodical for some $n\in{\mathbb N}$, then $\varphi_4-\varphi_3+\varphi_2-\varphi_1\in\pi{\mathbb Q},$ which implies the conditions (\ref{eq:30}) are satisfied.

Such a way we arrive to the following statement.

{\bf Theorem 5.} The periodicity of characteristic billiard on the unit disk is necessary and sufficient for the violation of Fredholm property of the Dirichlet problem $(\ref{eq:01})^0$, $(\ref{eq:28})^0$ in $L^2(K),$ and its kernel consists of countable set of linearly independent polynomial solutions (\ref{eq:30}).

\section*{Funding}
This work is supported by the Volkswagen Foundation (the project numbers are A131968 and 9C624) and by the Ministry of Education and Science of Ukraine (project number is 0121U109525).

\end{document}